\documentclass[12pt,twoside]{article}
\usepackage{amsmath, amssymb, amsthm}
\font\fourteenb=cmb10 at 14pt
\begin{document}
\begin{center}

\vspace*{-1.0cm}\noindent \
 \textbf{MATHEMATICS AND EDUCATION IN MATHEMATICS, 1986}\\[-0.0mm]\
\ Proceedings of 15th Spring Conference of

the Union of Bulgarian Mathematicians,\\[-0.0mm]
\textit{ Sunny Beach, April   2 - 6, 1986}\\[-0.0mm]
\font\fourteenb=cmb10 at 12pt \ \

{\bf \LARGE Multipliers of integrals of Cauchy - Stieltjes type
   \\ \ \\ \large Peyo Stoilov}
\end{center}

\

\footnotetext{

Theorem of Smirnov - Kotchine, Closed Jordan curve.}
\begin{abstract}
 Let  ${\rm {\mathbb G}}$  be a domain with closed rectifiable Jordan curve  $\ell $ .
 Let  $K({\rm {\mathbb G}})$  be the space of all analytic functions in  ${\rm {\mathbb G}}
 $  representable by a Cauchy - Stieltjes integral.
 Let  ${\rm {\mathfrak M}}(K)$  be the class of all multipliers of the space  $K({\rm {\mathbb G}}).$
 In this paper we prove that if  $f$  is bounded analytic function on  ${\rm {\mathbb G}}$  and
     $${\kern 1pt} {\kern 1pt} {\kern 1pt} \mathop{ess\sup }\limits_{\eta \in \ell } \int _{\ell }
     \frac{\left|f(\zeta )-f(\eta )\right|}{\left|\zeta -\eta \right|} \left|d\zeta \right|{\kern 1pt}
      {\kern 1pt} {\kern 1pt} {\kern 1pt} {\kern 1pt} {\kern 1pt} {\kern 1pt}
      {\kern 1pt} <\infty {\kern 1pt} {\kern 1pt} {\kern 1pt} {\kern 1pt} ,$$
then  $f\in {\rm {\mathfrak M}}(K)$ .

     If  ${\rm {\mathbb G}}={\rm {\mathbb D}}$  is the unit disc, this theorem was proved
for the first time by V. P. Havin. In particular for a smooth
curve $\ell $ we prove that if $f'\in E^{p} ({\rm {\mathbb
G}}),{\kern 1pt} {\kern 1pt} {\kern 1pt} {\kern 1pt} {\kern 1pt}
{\kern 1pt} p>1,$ then  $f\in {\rm {\mathfrak M}}(K),$  where
$E^{p} ({\rm {\mathbb G}})$  are the spaces of Smirnov.

\end{abstract}
\section{Introduction}

\

Let  $\ell $  be a closed rectifiable Jordan curve bounding the
finite and infinite domains  ${\rm {\mathbb G}}$  and  ${\rm
{\mathbb G}}^{-} $ {\it  }in the complex plane  ${\mathbb C} $ .
Let  $\ell _{r} $   $(0<r<1)$  be the map of circle  ${\kern 1pt}
\left|z\right|{\kern 1pt} =r$ {\it  }under conformal mapping
$\zeta =\varphi (z)$  of the unit disc  ${\rm {\mathbb D}}$  on
the domain  ${\rm {\mathbb G}}$ {\it . }Let denote by{\it  }
$A({\rm {\mathbb Q} })$ {\it  }the set of all analytic functions
in a domain  ${\rm {\mathbb Q} }\subset {\mathbb C} $ . For
${\kern 1pt} {\kern 1pt} {\kern 1pt} {\kern 1pt} {\kern 1pt}
0<p\le \infty $ , let  $E^{p} ({\rm {\mathbb G}})$  be the Smirnov
class [1,2] of functions  $f\in A({\rm {\mathbb G}})$ , for which

$$\left\| f\right\| _{E^{p} }^{p} =\mathop{\sup }\limits_{0<{\kern 1pt} {\kern 1pt} {\kern 1pt} r<{\kern 1pt} {\kern 1pt} {\kern 1pt} 1} {\kern 1pt} {\kern 1pt} {\kern 1pt} {\kern 1pt} {\kern 1pt} \frac{1}{2\pi } \int _{\ell _{r} }\left|f(\zeta )\right|^{p} \left|d\zeta \right| <\infty ,{\kern 1pt} {\kern 1pt} {\kern 1pt} {\kern 1pt} {\kern 1pt} {\kern 1pt} {\kern 1pt} {\kern 1pt} {\kern 1pt} {\kern 1pt} {\kern 1pt} {\kern 1pt} {\kern 1pt} {\kern 1pt} {\kern 1pt} {\kern 1pt} {\kern 1pt} {\kern 1pt} {\kern 1pt} {\kern 1pt} {\kern 1pt} {\kern 1pt} {\kern 1pt} {\kern 1pt} 0<p<\infty ,$$

             $$\left\| f\right\| _{E^{\infty } } =\mathop{\sup }\limits_{\zeta \in {\rm {\mathbb G}}} {\kern 1pt} {\kern 1pt} {\kern 1pt} {\kern 1pt} \left|f(\zeta )\right|<\infty ,{\kern 1pt} {\kern 1pt} {\kern 1pt} {\kern 1pt} {\kern 1pt} {\kern 1pt} {\kern 1pt} {\kern 1pt} {\kern 1pt} {\kern 1pt} {\kern 1pt} {\kern 1pt} {\kern 1pt} {\kern 1pt} {\kern 1pt} {\kern 1pt} {\kern 1pt} {\kern 1pt} {\kern 1pt} {\kern 1pt} {\kern 1pt} {\kern 1pt} {\kern 1pt} {\kern 1pt} p=\infty .$$

If   ${\rm {\mathbb G}}={\rm {\mathbb D}}$ ,  ${\rm {\mathbb
D}}=\left\{\zeta :{\kern 1pt} {\kern 1pt} {\kern 1pt} {\kern 1pt}
{\kern 1pt} {\kern 1pt} \left|\zeta \right|{\kern 1pt} <1\right\}$
,  then   $E^{p} ({\rm {\mathbb G}})=H^{p} $   is the usual Hardy
class.

\

Let  $N$ denote the Nevanlinna class:

\

     $f\in N$     $\Leftrightarrow $     $f\in A({\rm {\mathbb D}}):{\kern 1pt} {\kern 1pt} {\kern 1pt} {\kern 1pt} {\kern 1pt} $  $\mathop{\lim }\limits_{{\kern 1pt} {\kern 1pt} {\kern 1pt} {\kern 1pt} {\kern 1pt} {\kern 1pt} {\kern 1pt} {\kern 1pt} r\to 1^{-} } {\kern 1pt} {\kern 1pt} {\kern 1pt} \int \limits _{\left|\zeta \right|=1}\log ^{+} \left|f(r\zeta )\right| \left|d\zeta \right|<\infty $ ,

\

    and  $N^{+} $ - the Smirnov  class:

\

     $f\in N^{+} $     $\Leftrightarrow $     $f\in N:{\kern 1pt} {\kern 1pt} {\kern 1pt} $  $\mathop{{\kern 1pt} {\kern 1pt} {\kern 1pt} {\kern 1pt} {\kern 1pt} {\kern 1pt} {\kern 1pt} {\kern 1pt} {\kern 1pt} {\kern 1pt} \lim }\limits_{{\kern 1pt} {\kern 1pt} {\kern 1pt} {\kern 1pt} {\kern 1pt} {\kern 1pt} {\kern 1pt} {\kern 1pt} r\to 1^{-} } {\kern 1pt} {\kern 1pt} {\kern 1pt} \int \limits _{\left|\zeta \right|=1}\log ^{+} \left|f(r\zeta )\right| \left|d\zeta \right|=\int \limits _{\left|\zeta \right|=1}\log ^{+} \left|f(\zeta )\right| \left|d\zeta \right|$ .

\

     {\bf Some basic facts.}

     \

    {\bf 1.}  $f(\zeta )\in E^{p} ({\rm {\mathbb G}})$     $\Leftrightarrow $     $f(\varphi (z))\sqrt[{p}]{\varphi '(z)} \in H^{p} $ ;

\

     {\bf 2.} {\bf  }Theorem of Smirnov - Kotchine:

\

     {\it  If}{\it   } $f\in N^{+} $ {\it  and has } $L^{p} $ {\it  boundary values, then } $f\in H^{p} $ {\it .}

\

     Let  $M(\ell )$  is the space of all finite, complex Borel measures on  $\ell $  with the usual variation norm. We denote by  $K({\rm {\mathbb G}})$  the space of all functions  $g\in A({\rm {\mathbb G}}),$  representable by an integral of Cauchy - Stieltjes type

$$g(z)=\int _{\ell }{\kern 1pt} \frac{d\mu (\zeta )}{\zeta -z}  {\kern 1pt} {\kern 1pt} {\kern 1pt} {\kern 1pt} {\kern 1pt} {\kern 1pt} {\kern 1pt} {\kern 1pt} \mathop{=}\limits^{def} {\kern 1pt} {\kern 1pt} {\kern 1pt} {\kern 1pt} {\kern 1pt} {\kern 1pt} {\kern 1pt} {\kern 1pt} K_{\mu } (z) ,{\kern 1pt} {\kern 1pt} {\kern 1pt}{\kern 1pt} {\kern 1pt} {\kern 1pt}\mu \in M(\ell ).$$

     If  $g\in K({\rm {\mathbb G}})$ , then  $M_{g} $  will denote the set of all measures  $\mu \in M(\ell )$  participating in the representation of   $g$ .

     Since   $K_{\mu } (z)=0,{\kern 1pt} {\kern 1pt} {\kern 1pt} {\kern 1pt} {\kern 1pt} {\kern 1pt} \forall z\in {\rm {\mathbb G}}$    $\Leftrightarrow $  $K_{\mu } (z)\in E^{1} ({\rm {\mathbb G}}^{-} )$  and
$$\displaystyle E^{1} ({\rm {\mathbb G}}^{-} )=\left\{K_{\mu }
:{\kern 1pt} {\kern 1pt} {\kern 1pt} {\kern 1pt} {\kern 1pt} d\mu
=f(\zeta )d\zeta ,{\kern 1pt} {\kern 1pt} {\kern 1pt} {\kern 1pt}
{\kern 1pt} f(\zeta )\in L^{1} (\ell ){\kern 1pt} ,{\kern 1pt}
{\kern 1pt} {\kern 1pt} {\kern 1pt} \int _{\ell }\zeta ^{n}
f(\zeta ) {\kern 1pt} d\zeta =0,{\kern 1pt} {\kern 1pt} {\kern
1pt} {\kern 1pt} {\kern 1pt} n=-1,{\kern 1pt} {\kern 1pt} -2{\kern
1pt} {\kern 1pt} {\kern 1pt} ....{\kern 1pt} \right\},$$
then the
space   $K({\rm {\mathbb G}})$ {\it   }is the Banach space with
the natural norm
$$\displaystyle \left\| g\right\| _{K({\rm {\mathbb G}})} =\inf
\left\{\left\| \mu \right\| _{M(\ell )} :{\kern 1pt} {\kern 1pt}
{\kern 1pt} {\kern 1pt} {\kern 1pt} {\kern 1pt} {\kern 1pt} {\kern
1pt} {\kern 1pt} {\kern 1pt} {\kern 1pt} \mu \in M_{g} \right\},$$
isometrically  isomorphic  to  $M(\ell )/E^{1} ({\rm {\mathbb
G}}^{-} ).$

     Let  $C_{A} ({\rm {\mathbb G}}^{-} )$ denote the space of all  functions   $h$ , analytic in  ${\rm {\mathbb G}}^{-} $  and continuous on  $\overline{{\rm {\mathbb G}}^{-} }={\rm {\mathbb G}}^{-} \cup \ell $ , with norm

$$\left\| h\right\| _{\infty } =\sup {\kern 1pt} {\kern 1pt} {\kern 1pt} \left\{{\kern 1pt} \left|h(z)\right|:{\kern 1pt} {\kern 1pt} {\kern 1pt} {\kern 1pt} {\kern 1pt} {\kern 1pt} {\kern 1pt} {\kern 1pt} {\kern 1pt} {\kern 1pt} {\kern 1pt} z\in {\rm {\mathbb G}}^{-} \right\}.$$

    Since   $M(\ell )/E^{1} ({\rm {\mathbb G}}^{-} )\cong C_{A} ^{*} ({\rm {\mathbb G}}^{-} )$ ,  then  $K({\rm {\mathbb G}})\cong C_{A} ^{*} ({\rm {\mathbb G}}^{-} )$   and
     $$\left\| g\right\| _{K({\rm {\mathbb G}})} =\sup \left\{\left|\int _{\ell }h(\eta )d\mu (\eta ) \right|:{\kern 1pt} {\kern 1pt} {\kern 1pt} {\kern 1pt} {\kern 1pt} {\kern 1pt} {\kern 1pt} {\kern 1pt} {\kern 1pt} {\kern 1pt} {\kern 1pt} h\in C_{A} ({\rm {\mathbb G}}^{-} ),{\kern 1pt} {\kern 1pt} {\kern 1pt} {\kern 1pt} {\kern 1pt} \left\| h\right\| _{\infty } \le 1\right\} ,
      {\kern 1pt} {\kern 1pt} {\kern 1pt}{\kern 1pt} {\kern 1pt} {\kern 1pt}       \mu \in M_{g} .$$

      Taking into
       account that the set  $A_{0} =A(\overline{{\rm {\mathbb G}}^{-} }{\kern 1pt} )$  is dense in   $C_{A} ({\rm {\mathbb G}}^{-} )$ , then easily we receive for   $h\in A_{0} ,{\kern 1pt} {\kern 1pt} {\kern 1pt} {\kern 1pt} g=K_{\mu } $ :

     $$\left|\int _{\ell _{r} }h(\zeta )g(\zeta )d\zeta  \right|{\kern 1pt} {\kern 1pt} {\kern 1pt} {\kern 1pt} ={\kern 1pt} {\kern 1pt} {\kern 1pt} \left|\int _{\ell _{r} }h(\zeta )\left(\int _{\ell }{\kern 1pt} \frac{d\mu (\eta )}{\eta -\zeta }  \right)d\zeta  \right|=$$

     \

     $${\kern 1pt} {\kern 1pt} ={\kern 1pt} {\kern 1pt} {\kern 1pt} \left|\int _{\ell }{\kern 1pt}  \left(\int _{\ell _{r} }\frac{h(\zeta )}{\zeta -\eta } d\zeta  \right)d\mu (\eta )\right|=\left|\int _{\ell }{\kern 1pt}  h(\eta )d\mu (\eta )\right|$$
and consequently
    $$\left\| g\right\| _{K({\rm {\mathbb G}})} =\sup \left\{\mathop{\lim }\limits_{r\to 1^{-} } \left|\int _{\ell _{r} }h(\zeta )g(\zeta )d\zeta  \right|:{\kern 1pt} {\kern 1pt} {\kern 1pt} {\kern 1pt} {\kern 1pt} {\kern 1pt} {\kern 1pt} {\kern 1pt} {\kern 1pt} {\kern 1pt} {\kern 1pt} h\in A_{0} ,{\kern 1pt} {\kern 1pt} {\kern 1pt} {\kern 1pt} {\kern 1pt} \left\| h\right\| _{\infty } \le 1\right\}.$$

     Let  ${\rm {\mathfrak M}}(K)$  be the class to all multipliers of the space $K({\rm {\mathbb G}})$ :

     $${\rm {\mathfrak M}}(K)=\left\{f\in A({\rm {\mathbb G}}):{\kern 1pt} {\kern 1pt} {\kern 1pt} {\kern 1pt} {\kern 1pt} {\kern 1pt} {\kern 1pt} {\kern 1pt} {\kern 1pt} f{\kern 1pt} g\in K({\rm {\mathbb G}}),{\kern 1pt} {\kern 1pt} {\kern 1pt} {\kern 1pt} {\kern 1pt} {\kern 1pt} {\kern 1pt} {\kern 1pt} {\kern 1pt} {\kern 1pt} \forall {\kern 1pt} g\in K({\rm {\mathbb G}})\right\}.$$

     \

     {\bf Lemma. }{\it If } $f\in {\rm {\mathfrak M}}(K)$ {\it , then}

     $$\sup \left\{{\kern 1pt} {\kern 1pt} {\kern 1pt} {\kern 1pt} {\kern 1pt} \left\| f{\kern 1pt} g\right\| _{K({\rm {\mathbb G}})} {\kern 1pt} {\kern 1pt} {\kern 1pt} :{\kern 1pt} {\kern 1pt} {\kern 1pt} {\kern 1pt} {\kern 1pt} {\kern 1pt} {\kern 1pt} g{\kern 1pt} \in K({\rm {\mathbb G}}),{\kern 1pt} {\kern 1pt} {\kern 1pt} {\kern 1pt} {\kern 1pt} {\kern 1pt} \left\| {\kern 1pt} g\right\| _{K({\rm {\mathbb G}})} \le 1{\kern 1pt} {\kern 1pt} {\kern 1pt} {\kern 1pt} {\kern 1pt} {\kern 1pt} \right\}<\infty .$$

     {\it Proof.} Let $f\in {\rm {\mathfrak M}}(K).$  We will show that the operator    $g{\kern 1pt} {\kern 1pt} {\kern 1pt} {\kern 1pt} {\kern 1pt} {\kern 1pt} {\kern 1pt} \to {\kern 1pt} {\kern 1pt} {\kern 1pt} {\kern 1pt} {\kern 1pt} {\kern 1pt} {\kern 1pt} fg$    is bounded in  $K({\rm {\mathbb G}}).$  By the closed graph theorem it suffices to show that the graph of the operator   $g{\kern 1pt} {\kern 1pt} {\kern 1pt} {\kern 1pt} {\kern 1pt} {\kern 1pt} {\kern 1pt} \to {\kern 1pt} {\kern 1pt} {\kern 1pt} {\kern 1pt} {\kern 1pt} {\kern 1pt} {\kern 1pt} fg$    is closed on  $K({\rm {\mathbb G}}).$

     Suppose   $\left(g_{n} \right)_{n\ge 0} \subset K({\rm {\mathbb G}})$   with   $g_{n} {\kern 1pt} {\kern 1pt} {\kern 1pt} {\kern 1pt} {\kern 1pt} {\kern 1pt} {\kern 1pt} \to {\kern 1pt} {\kern 1pt} {\kern 1pt} {\kern 1pt} {\kern 1pt} {\kern 1pt} {\kern 1pt} g$   and   $fg_{n} {\kern 1pt} {\kern 1pt} {\kern 1pt} {\kern 1pt} {\kern 1pt} {\kern 1pt} {\kern 1pt} \to {\kern 1pt} {\kern 1pt} {\kern 1pt} {\kern 1pt} {\kern 1pt} {\kern 1pt} {\kern 1pt} \phi $   in the  $K({\rm {\mathbb G}})$  - norm.

     {\it If} $F\in K({\rm {\mathbb G}})$ {\it , }then for all $z\in {\rm {\mathbb G}}$ ,

      $$F(z)=\int _{\ell }{\kern 1pt} \frac{d\mu (\zeta )}{\zeta -z}{\kern 1pt}{\kern 1pt} {\kern 1pt} {\kern 1pt} {\kern 1pt}              \Rightarrow                {\kern 1pt} {\kern 1pt} {\kern 1pt} {\kern 1pt}\left|F(z)\right|\le \left(\inf \left|\zeta -z\right|:\zeta \in \ell \right)^{-1} \left\| {\kern 1pt} F\right\| _{K({\rm {\mathbb G}})} .$$

     From this inequality follows that the converse in norm in $K({\rm {\mathbb G}})$  implies a pointwise converse on ${\rm {\mathbb G}}$ . Consequently  $\phi =fg$  and the graph of the operator    $g{\kern 1pt} {\kern 1pt} {\kern 1pt} {\kern 1pt} {\kern 1pt} {\kern 1pt} {\kern 1pt} \to {\kern 1pt} {\kern 1pt} {\kern 1pt} {\kern 1pt} {\kern 1pt} {\kern 1pt} {\kern 1pt} fg$    is closed on  $K({\rm {\mathbb
     G}}).$ $\Box$

     From proved lemma follows that  ${\rm {\mathfrak M}}(K)$ is commutative Banach algebra with norm defined by

     $$\left\| f\right\| _{{\rm {\mathfrak M}}(K)} =\sup \left\{{\kern 1pt} {\kern 1pt} {\kern 1pt} {\kern 1pt} {\kern 1pt} \left\| f{\kern 1pt} g\right\| _{K({\rm {\mathbb G}})} {\kern 1pt} {\kern 1pt} {\kern 1pt} :{\kern 1pt} {\kern 1pt} {\kern 1pt} {\kern 1pt} {\kern 1pt} {\kern 1pt} {\kern 1pt} g{\kern 1pt} \in K({\rm {\mathbb G}}),{\kern 1pt} {\kern 1pt} {\kern 1pt} {\kern 1pt} {\kern 1pt} {\kern 1pt} \left\| {\kern 1pt} g\right\| _{K({\rm {\mathbb G}})} \le 1{\kern 1pt} {\kern 1pt} {\kern 1pt} {\kern 1pt} {\kern 1pt} {\kern 1pt} \right\}.$$

     {\it If } ${\rm {\mathbb G}}={\rm {\mathbb D}}$ {\it , }then the study of the space  ${\rm {\mathfrak M}}={\rm {\mathfrak M}}(K)$  was started in papers of V. P. Havin [3, 4] and was continued in [5, 6, 7]. For example, in [5, 6] was proved:

     {\it If } $f$ {\it  is bounded} {\it  analytic function on } ${\rm {\mathbb D}}$ {\it  and}

     $${\kern 1pt} {\kern 1pt} {\kern 1pt} \mathop{ess\sup }\limits_{\left|\eta \right|=1} \int _{\left|\zeta \right|=1} \frac{\left|f(\zeta )-f(\eta )\right|}{\left|\zeta -\eta \right|} \left|d\zeta \right|{\kern 1pt} {\kern 1pt} {\kern 1pt} {\kern 1pt} {\kern 1pt} {\kern 1pt} {\kern 1pt} {\kern 1pt} <\infty {\kern 1pt} {\kern 1pt} {\kern 1pt} {\kern 1pt} ,$$

     {\it then} $f\in {\rm {\mathfrak M}}$ {\it .}

\

     In this paper this result is generalized for the multipliers of the space $K({\rm {\mathbb G}}).$
\section{Main results}
 {\bf Theorem 1.} {\it Let } ${\rm {\mathbb G}}$ {\it  }{\it be a domain with closed rectifiable Jordan curve} $\ell $ {\it .}

     {\it If } $f\in E^{\infty } ({\rm {\mathbb G}})$ {\it  and}
     $${\kern 1pt} \Lambda (f){\kern 1pt} {\kern 1pt} {\kern 1pt} {\kern 1pt} {\kern 1pt} \mathop{=}\limits^{def} {\kern 1pt} {\kern 1pt} {\kern 1pt} {\kern 1pt} {\kern 1pt} {\kern 1pt} {\kern 1pt} \mathop{ess\sup }\limits_{\eta \in \ell } \int _{\ell } \frac{\left|f(\zeta )-f(\eta )\right|}{\left|\zeta -\eta \right|} \left|d\zeta \right|{\kern 1pt} {\kern 1pt} {\kern 1pt} {\kern 1pt} {\kern 1pt} {\kern 1pt} {\kern 1pt} {\kern 1pt} <\infty {\kern 1pt} {\kern 1pt} {\kern 1pt} {\kern 1pt} ,$$
{\it  then  } $f\in {\rm {\mathfrak M}}(K)$ {\it   and   }
$$\left\| f\right\| _{{\rm {\mathfrak M}}(K)} \le \left\| f\right\|
_{E^{\infty } } +\Lambda (f).$$

\

     {\it Proof. }Let  $f\in E^{\infty } ({\rm {\mathbb G}})$ {\it  }and $\Lambda (f){\kern 1pt} {\kern 1pt} {\kern 1pt} {\kern 1pt}
 <\infty $ . Let  ${\rm {\mathbb E}}\subseteq \ell $  be a subset with total Lebesgue measure  $(m({\rm {\mathbb E}})=m(\ell ))$  lying on  $\ell $  such that
$$\left\| f\right\| _{E^{\infty } } =\mathop{\sup }\limits_{\eta \in {\rm {\mathbb E}}{\kern 1pt} } {\kern 1pt} {\kern 1pt} {\kern 1pt} \left|f(\eta )\right|.$$
Then
$$\left\| f\right\| _{{\rm {\mathfrak M}}(K)} =\sup \left\{{\kern 1pt} {\kern 1pt} {\kern 1pt} {\kern 1pt} {\kern 1pt}
\left\| f{\kern 1pt} g\right\| _{K({\rm {\mathbb G}})} {\kern 1pt} {\kern 1pt} {\kern 1pt} :{\kern 1pt} {\kern 1pt} {\kern 1pt} {\kern 1pt} {\kern 1pt} {\kern 1pt} {\kern 1pt} g{\kern 1pt} \in K({\rm {\mathbb G}}),{\kern 1pt} {\kern 1pt} {\kern 1pt} {\kern 1pt} {\kern 1pt} {\kern 1pt} \left\| {\kern 1pt} g\right\| _{K({\rm {\mathbb G}})} \le 1{\kern 1pt} {\kern 1pt} {\kern 1pt} {\kern 1pt} {\kern 1pt} {\kern 1pt} \right\}=$$
$$\displaystyle =\sup \left\{\mathop{\lim }\limits_{r\to 1^{-} }
\left|\int _{\ell _{r} }f(\zeta )g(\zeta )h(\zeta )d\zeta
\right|:{\kern 1pt} {\kern 1pt} {\kern 1pt} {\kern 1pt} {\kern
1pt} {\kern 1pt} {\kern 1pt} g{\kern 1pt} \in K({\rm {\mathbb
G}}),{\kern 1pt} {\kern 1pt} {\kern 1pt} {\kern 1pt} {\kern 1pt}
{\kern 1pt} \left\| {\kern 1pt} g\right\| _{K({\rm {\mathbb G}})}
\le 1,{\kern 1pt} {\kern 1pt} {\kern 1pt} {\kern 1pt} {\kern 1pt}
{\kern 1pt} {\kern 1pt} {\kern 1pt} h\in A_{0} ,{\kern 1pt} {\kern
1pt} {\kern 1pt} {\kern 1pt} {\kern 1pt} \left\| h\right\|
_{\infty } \le 1\right\},$$ where  $\ell _{r} $   $(0<r<1)$   is
the map of circle  ${\kern 1pt} \left|z\right|{\kern 1pt} =r$ {\it
}under conformal mapping $\zeta =\varphi (z)$  of the unit disc
${\rm {\mathbb D}}$  on the domain  ${\rm {\mathbb G}}.$

     If $g\in K({\rm {\mathbb G}})$ ,  then
      $$g(z)=\int _{\ell }{\kern 1pt} \frac{d\mu (\eta )}{\eta -z},              z\in {\rm {\mathbb G}}$$
     and
     $$\left|\int _{\ell _{r} }f(\zeta )g(\zeta )h(\zeta )d\zeta  \right|{\kern 1pt} {\kern 1pt} {\kern 1pt} {\kern 1pt} ={\kern 1pt} {\kern 1pt} {\kern 1pt} \left|\int _{\ell _{r} }f(\zeta )h(\zeta )\left(\int _{\ell }{\kern 1pt} \frac{d\mu (\eta )}{\eta -\zeta }  \right)d\zeta  \right|=$$

\

     $${\kern 1pt} {\kern 1pt} ={\kern 1pt} {\kern 1pt} {\kern 1pt} \left|\int _{\ell }{\kern 1pt}  \left(\int _{\ell _{r} }\frac{f(\zeta )h(\zeta )}{\zeta -\eta } d\zeta  \right)d\mu (\eta )\right|\le $$

\

     $$\le \int _{\ell }{\kern 1pt} \left|d\mu (\eta )\right| {\kern 1pt} {\kern 1pt} {\kern 1pt} {\kern 1pt} \sup {\kern 1pt} {\kern 1pt} \left\{{\kern 1pt} \left|\int _{\ell _{r} }\frac{f(\zeta )h(\zeta )}{\zeta -\eta } d\zeta  \right|:{\kern 1pt} {\kern 1pt} {\kern 1pt} {\kern 1pt} {\kern 1pt} {\kern 1pt} {\kern 1pt} {\kern 1pt} \eta \in {\rm {\mathbb E}}\right\}{\kern 1pt} {\kern 1pt}, $$
    $\Rightarrow $
         $$\left|\int _{\ell _{r} }f(\zeta )g(\zeta )h(\zeta )d\zeta  \right|{\kern 1pt} {\kern 1pt} {\kern 1pt} \le {\kern 1pt} {\kern 1pt} {\kern 1pt} \left\| {\kern 1pt} g\right\| _{K({\rm {\mathbb G}})} \sup {\kern 1pt} {\kern 1pt} \left\{{\kern 1pt} \left|\int _{\ell _{r} }\frac{f(\zeta )h(\zeta )}{\zeta -\eta } d\zeta  \right|:{\kern 1pt} {\kern 1pt} {\kern 1pt} {\kern 1pt} {\kern 1pt} {\kern 1pt} {\kern 1pt} {\kern 1pt} \eta \in {\rm {\mathbb E}}\right\}.$$

     Then
     $$\left\| f\right\| _{{\rm {\mathfrak M}}(K)} \le \sup {\kern 1pt} {\kern 1pt} \left\{\mathop{\lim }\limits_{r\to 1^{-} } {\kern 1pt} \left|\int _{\ell _{r} }\frac{f(\zeta )h(\zeta )}{\zeta -\eta } d\zeta  \right|:{\kern 1pt} {\kern 1pt} {\kern 1pt} {\kern 1pt} {\kern 1pt} {\kern 1pt} {\kern 1pt} {\kern 1pt} \eta \in {\rm {\mathbb E}},{\kern 1pt} {\kern 1pt} {\kern 1pt} {\kern 1pt} {\kern 1pt} {\kern 1pt} h\in A_{0} ,{\kern 1pt} {\kern 1pt} {\kern 1pt} {\kern 1pt} {\kern 1pt} \left\| h\right\| _{\infty } \le 1{\kern 1pt} \right\}.$$

     Further, if  $\eta \in {\rm {\mathbb E}},{\kern 1pt} {\kern 1pt} {\kern 1pt} {\kern 1pt} {\kern 1pt} {\kern 1pt} h\in A_{0} ,{\kern 1pt} {\kern 1pt} {\kern 1pt} {\kern 1pt} {\kern 1pt} \left\| h\right\| _{\infty } \le 1$ , we have

     $${\kern 1pt} \left|\int _{\ell _{r} }\frac{f(\zeta )h(\zeta )}{\zeta -\eta } d\zeta  \right|\le \left|\int _{\ell _{r} }\frac{f(\zeta )-f(\eta )}{\zeta -\eta } h(\zeta ){\kern 1pt} {\kern 1pt} d\zeta  \right|+\left|f(\eta )\right|\left|\int _{\ell _{r} }\frac{h(\zeta )}{\zeta -\eta } d\zeta  \right|\le $$

     $${\kern 1pt} \le \left\| h\right\| _{L^{\infty } (\ell _{r} )} \int _{\ell _{r} }\left|\frac{f(\zeta )-f(\eta )}{\zeta -\eta } \right|\left|d\zeta \right| +\left\| f\right\| _{E^{\infty } } .$$

     Then
     $$\left\| f\right\| _{{\rm {\mathfrak M}}(K)} \le \sup {\kern 1pt} {\kern 1pt} \left\{\mathop{\lim }\limits_{r\to 1^{-} } {\kern 1pt} \int _{\ell _{r} }\left|\frac{f(\zeta )-f(\eta )}{\zeta -\eta } \right|\left|d\zeta \right| :{\kern 1pt} {\kern 1pt} {\kern 1pt} {\kern 1pt} {\kern 1pt} {\kern 1pt} {\kern 1pt} {\kern 1pt} \eta \in {\rm {\mathbb E}}{\kern 1pt} {\kern 1pt} {\kern 1pt} {\kern 1pt} {\kern 1pt} {\kern 1pt} {\kern 1pt} \right\}+\left\| f\right\| _{E^{\infty } } .$$

     We denote for  $\eta \in {\rm {\mathbb E}}$
       $$F_{\eta } (\zeta )=\frac{f(\zeta )-f(\eta )}{\zeta -\eta } ,{\kern 1pt} {\kern 1pt} {\kern 1pt} {\kern 1pt} {\kern 1pt} {\kern 1pt} {\kern 1pt} {\kern 1pt} {\kern 1pt} {\kern 1pt} \zeta \in {\rm {\mathbb G}}.$$

     Then
     $$\left\| f\right\| _{{\rm {\mathfrak M}}(K)} \le \sup {\kern 1pt} {\kern 1pt} \left\{\left\| F_{\eta } (\zeta )\right\| _{E^{1} } :{\kern 1pt} {\kern 1pt} {\kern 1pt} {\kern 1pt} {\kern 1pt} {\kern 1pt} {\kern 1pt} {\kern 1pt} \eta \in {\rm {\mathbb E}}{\kern 1pt} {\kern 1pt} {\kern 1pt} {\kern 1pt} {\kern 1pt} {\kern 1pt} {\kern 1pt} \right\}+\left\| f\right\| _{E^{\infty } } .$$

\

     To end the proof it is necessary to show that
     $$\Lambda (f){\kern 1pt} {\kern 1pt} {\kern 1pt} {\kern 1pt} <\infty {\kern 1pt} {\kern 1pt} {\kern 1pt} {\kern 1pt} {\kern 1pt} {\kern 1pt} {\kern 1pt} {\kern 1pt} \Rightarrow {\kern 1pt} {\kern 1pt} {\kern 1pt} {\kern 1pt} {\kern 1pt} \mathop{\sup }\limits_{\eta \in {\rm {\mathbb E}}} \left\| F_{\eta } \right\| _{E^{1} } <\infty .{\kern 1pt} $$

     Let  $\zeta =\varphi (z)$  is conformal mapping of the unit disc  ${\rm {\mathbb D}}$  on the domain ${\rm {\mathbb G}}$ ,  $\eta \in {\rm {\mathbb E}}$
and  $\eta =\varphi (t)$ ,   $\left|t\right|=1$ . Then

      $$F_{\eta } (\zeta )=\frac{f(\zeta )-f(\eta )}{\zeta -\eta } \in E^{1} ({\rm {\mathbb G}}){\kern 1pt} {\kern 1pt} {\kern 1pt}        \Leftrightarrow        {\kern 1pt} {\kern 1pt} {\kern 1pt}\Omega (z){\kern 1pt} {\kern 1pt} {\kern 1pt} {\kern 1pt} {\kern 1pt} \mathop{=}\limits^{def} {\kern 1pt} {\kern 1pt} {\kern 1pt} {\kern 1pt} {\kern 1pt} {\kern 1pt} \frac{f(\varphi (z))-f(\varphi (t))}{\varphi (z)-\varphi (t)} \varphi '(z)\in H^{1} .$$

     In [8] was proved that the functions
     $$\displaystyle \frac{1}{\varphi (z)-\varphi (t)} \in N^{+} .$$

     Since
      $$f(\varphi (z))-f(\varphi (t))\in H^{\infty },      {\kern 1pt} {\kern 1pt} \varphi '(z)\in H^{1} ,{\kern 1pt} {\kern 1pt} {\kern 1pt} {\kern 1pt} {\kern 1pt} {\kern 1pt}     then {\kern 1pt} {\kern 1pt} {\kern 1pt} {\kern 1pt} {\kern 1pt} {\kern 1pt} \Omega (z){\kern 1pt} {\kern 1pt} {\kern 1pt} {\kern 1pt} {\kern 1pt} \in N^{+} .$$

     Besides
      $$\int _{\left|z\right|=1}\left|\Omega (z)\right| \left|dz\right| =  \int _{\ell }\left|F_{\eta } (\zeta )\right| \left|d\zeta \right|\le \Lambda (f){\kern 1pt} {\kern 1pt} {\kern 1pt} {\kern 1pt} <\infty {\kern 1pt} {\kern 1pt} $$
and according to the Theorem of Smirnov - Kotchine

 $$\Omega (z){\kern 1pt} {\kern 1pt} {\kern 1pt} {\kern 1pt} {\kern 1pt} \in H^{1} ,{\kern 1pt} {\kern 1pt} {\kern 1pt} {\kern 1pt} {\kern 1pt}   \left\| \Omega (z)\right\| _{H^{1} } \le \Lambda (f)<\infty .$$

     Consequently
      $$\left\| F_{\eta } \right\| _{E^{1} } =\left\| \Omega (z)\right\| _{H^{1} } \le \Lambda (f)<\infty {\kern 1pt} $$
      and
     $$\left\| f\right\| _{{\rm {\mathfrak M}}(K)} \le \Lambda (f)+\left\| f\right\| _{E^{\infty } } <\infty .$$$\Box$

     {\bf Theorem 2.} {\it Let } ${\rm {\mathbb G}}$ {\it  }{\it be a domain with a smooth curve } $\ell $ {\it .}

\

     {\it  If } $f\in E^{\infty } ({\rm {\mathbb G}})$ {\it  and  } $f'\in E^{p} ({\rm {\mathbb G}}),{\kern 1pt} {\kern 1pt} {\kern 1pt} {\kern 1pt} {\kern 1pt} {\kern 1pt} {\kern 1pt} {\kern 1pt} p>1,$ {\it  then } $f\in {\rm {\mathfrak M}}(K).$

\

     {\it Proof. }Let  $f'\in E^{p} ({\rm {\mathbb G}}),{\kern 1pt} {\kern 1pt} {\kern 1pt} {\kern 1pt} {\kern 1pt} {\kern 1pt} {\kern 1pt} {\kern 1pt} p>1,$ {\it  } $z=z(s){\kern 1pt} {\kern 1pt} {\kern 1pt} {\kern 1pt} {\kern 1pt} {\kern 1pt} {\kern 1pt} (0\le s\le s_{0} )$ {\it  }is the equation of  $\ell ,$ where the arc length is as parameter.
     We shall prove that

     $${\kern 1pt} \Lambda (f){\kern 1pt} {\kern 1pt} {\kern 1pt} {\kern 1pt} {\kern 1pt} ={\kern 1pt} {\kern 1pt} {\kern 1pt} {\kern 1pt} {\kern 1pt} {\kern 1pt} {\kern 1pt} \mathop{\sup }\limits_{\sigma \in \left[0,{\kern 1pt} {\kern 1pt} s_{0} \right]} {\kern 1pt} {\kern 1pt} {\kern 1pt} \int _{0}^{s_{0} }\frac{\left|f(s)-f(\sigma )\right|}{\left|z(s)-z(\sigma )\right|} ds{\kern 1pt} {\kern 1pt} {\kern 1pt} {\kern 1pt}  {\kern 1pt} <\infty {\kern 1pt} {\kern 1pt} {\kern 1pt} {\kern 1pt} ,$$
where  $f(s)=f\circ z(s).$

     Since  $\ell $  is a smooth curve, then
$$\left|z(s)-z(\sigma )\right|\ge \left|s-\sigma \right|,$$
where $c_{0} $  is a constant. Then

$${\kern 1pt} {\kern 1pt} \int _{0}^{s_{0} }\frac{\left|f(s)-f(\sigma )\right|}{\left|z(s)-z(\sigma )\right|} ds{\kern 1pt} {\kern 1pt} {\kern 1pt} {\kern 1pt}  \le \frac{1}{c_{0} } \int _{0}^{s_{0} }\frac{1}{\left|s-\sigma \right|} \left|\int _{\sigma }^{s}f'(x)dx{\kern 1pt} {\kern 1pt} {\kern 1pt} {\kern 1pt}  \right|ds{\kern 1pt} {\kern 1pt} {\kern 1pt} \le {\kern 1pt}  $$

$${\kern 1pt} {\kern 1pt} \le \frac{1}{c_{0} } \int _{0}^{s_{0} }\left(\int _{\sigma }^{s}\left|f'(x)\right|^{p} dx{\kern 1pt} {\kern 1pt} {\kern 1pt} {\kern 1pt}  \right)^{1/p} \left|s-\sigma \right|^{-1/p} ds{\kern 1pt} {\kern 1pt} {\kern 1pt} {\kern 1pt} {\kern 1pt} {\kern 1pt} {\kern 1pt} {\kern 1pt} {\kern 1pt} \le {\kern 1pt}  {\kern 1pt} {\kern 1pt} {\kern 1pt} {\kern 1pt} {\kern 1pt} Const.\left\| f'\right\| _{E^{p} ({\rm {\mathbb G}})} .$$

     Consequently
     $$\Lambda (f){\kern 1pt} {\kern 1pt} {\kern 1pt} {\kern 1pt} {\kern 1pt} \le {\kern 1pt} {\kern 1pt} {\kern 1pt} {\kern 1pt} {\kern 1pt} Const.\left\| f'\right\| _{E^{p} ({\rm {\mathbb G}})} {\kern 1pt} <\infty $$
and by Theorem 1  $f\in {\rm {\mathfrak M}}(K).$ $\Box$

\

    {\bf Remark. }\emph{It should be noted that if } ${\rm {\mathbb G}}={\rm {\mathbb D}}$ , \emph{Theorem 2 remains valid and for}  $p=1$
     \emph{( Theorem of S. A.Vinogradov }[5]:  $f'(z)\in H^{1} \Rightarrow $  $f\in {\rm {\mathfrak M}}$ \emph{)}.

   \emph{ Whether Theorem 2 for } $p=1$  \emph{is correct, if } ${\rm {\mathbb G}}$ \emph{ is domain with a smooth boundary, remains unknown.}

\

\

\noindent
{\small Department of Mathematics\\
        Technical University\\
        25, Tsanko Dustabanov,\\
        Plovdiv, Bulgaria\\
        e-mail: peyyyo@mail.bg}


\begin{thebibliography}{199}


\bibitem{1}I. I. Privalov. {\it Boundary properties of analytic
functions}. GATT, Moscow, 1950. (Russian)  German transl., VEB
Deutscher Verlag, Berlin, 1956. MR0047765 (13, 926h)

\bibitem{2}P. L. Duren. {\it Theory of } $H^{p} $ {\it  spaces}. Academic
Press, New York, 1970. MR0268655 (42 \#3552)

\bibitem{3} V. P. Havin. {\it On analytic functions representable by an
integral of Cauchy-Stieltjes type}. Vestnik Leningrad. Univ. Ser.
Mat. Meh. Astronom. 13 no.1, 1958, 66-79. (Russian) MR0095256 (20
\#1762)

\bibitem{4}V. P. Havin. {\it Relations between certain classes of functions
regular in the unit disk.} Vestnik Leningrad. Univ. Ser. Mat. Meh.
Astronom. 17 no.1, 1962, 102-110. (Russian)  MR0152660 (27 \#2635)

\bibitem{5}S. A. Vinogradov. {\it Properties of multipliers of Cauchy -
Stieltjes integrals and some factorization problems for analytic
functions.} Amer. Math. Sos. Transl. (2) vol. 115, 1980, 1-32.
MR0586560 (58 \#28518)

\bibitem{6} S. A. Vinogradov, M. G. Goluzina, V. P. Havin. {\it
Multipliers}{\it  and divisors of Cauchy-Stieltjes type
integrals.} Zap. Nauen. Sem. Leningrad. Otdel. Mat. Inst. Steklov.
(LOMI) 19, 1970, 55-78, 1970. (Russian) MR0291471 (45 \#562)

\bibitem{7} S. V. Hru\v{s}\v{c}\"{e}v, S. A. Vinogradov. {\it Inner functions and
multipliers of Cauchy type integrals}. Ark. mat, 19, 1981, 23-42.
MR0625535 (83c:30027)

\bibitem{8}G. C. Tumarkin. {\it Properties of analytic functions representable by
integrals of Cauchy-Stieltjes and Cauchy-Lebesgue type.} Izv.
Akad. Nauk Armjan. SSR Ser. Fiz.- Mat. Nauk  16{\bf  } no. 5,
1963, 23-45. (Russian) MR0156955 (28 \#197)

\end{thebibliography}
\end{document}